\RequirePackage{lineno}
\documentclass[11pt]{article}

\usepackage{amsmath,amsfonts,amsbsy,amstext,amsopn,amssymb}
\usepackage{longtable}
\usepackage{theorem}
\usepackage{amsfonts,amsbsy,amstext,amsopn}
\usepackage{graphicx}
\usepackage{anysize,indentfirst,geometry,xcolor}
\usepackage[colorlinks,citecolor=cyan,linkcolor=blue]{hyperref}
\usepackage{cite}

\makeatletter \@addtoreset{equation}{section}

\makeatother

\usepackage{latexsym}
\usepackage{amssymb}
\usepackage{amsmath}

\usepackage{amsmath,amsfonts,amsbsy,amstext,amsopn,amssymb}
\usepackage{longtable}
\usepackage{theorem}

\newtheorem{theorem}{Theorem}

\newtheorem{lemma}{Lemma}
\newtheorem{corollary}{Corollary}

\newcommand{\R}{{\mathbb{R}}}
\newcommand{\C}{{\mathbb{C}}}
\newcommand{\HM}{{\mathbb{H}}}
\newcommand{\SM}{{\mathbb{S}}}
\newcommand{\SK}{{\mathbb{SK}}}
\newcommand{\sym}{{\sf{sym}}}
\newcommand{\ske}{{\sf{skew}}}

\newcommand{\diag}{{\sf{Diag}}}

\renewcommand{\Im}{{\rm{Im}}}

\def\e{\varepsilon}
\def\Prob{\mathop{\rm Prob}}

\newcommand{\bfi}{{\mathbf i}}

\newcommand{\cals}{{\mathcal S}}
\newcommand{\calm}{{\mathcal M}}
\newcommand{\caln}{{\mathcal N}}
\newcommand{\cald}{{\mathcal D}}
\newcommand{\bfe}{{\bf e}}

\newcommand\qed{{\unskip\nobreak\hfil\penalty50\hskip2em\vadjust{}
\nobreak\hfil$\Box$\parfillskip=0pt\finalhyphendemerits=0\par}}

\def\Oh{{\mathcal O}}

\newcommand{\vect}{{\sf{vec}}}
\def\Re{{\sf Re}}
\def\Im{{\sf Im}}

\begin{document}


\title{\bf  Structured condition numbers and small sample
condition estimation of symmetric algebraic Riccati
equations}
\author{
Huai-An Diao\thanks{ School of Mathematics and Statistics, Northeast Normal
University, Chang Chun 130024, P.R. China. Corresponding author. ({\tt hadiao@nenu.edu.cn and hadiao78@yahoo.com}) }
\and
Dongmei Liu\thanks{School of Mathematics and Statistics, Northeast Normal University,
Chang Chun 130024, P.R. China. ({\tt 627806414@qq.com}) Current address: No.1 Senior Middle School Of Inner Mongolia First Machinery  Group Corporation, Wenhua Rd. No. 26,  Baotou 014030, P.R. China}
\and
Sanzheng Qiao \thanks{Department of Computing and Software,
McMaster University, Hamilton, Ontario, L8S 4K1 Canada. ({\tt qiao@mcmaster.ca}) }
}

\date{}
\maketitle

\begin{quote}
{\bf Abstract.}  This paper is devoted to a structured perturbation analysis of
the symmetric algebraic Riccati equations by exploiting the
symmetry structure. Based on the analysis, the
upper bounds for the structured normwise, mixed and
componentwise condition numbers are derived.
Due to the exploitation of the symmetry structure, our results
are improvements of the previous work on the perturbation
analysis and condition numbers of the symmetric algebraic
Riccati equations. Our preliminary numerical experiments
demonstrate that our condition numbers provide accurate estimates
for the change in the solution caused by the perturbations
on the data.
Moreover, by applying the small sample condition estimation
method, we propose a statistical algorithm for practically estimating
the condition numbers of the symmetric algebraic Riccati
equations.
\end{quote}

{\small {\bf Keywords:} Symmetric algebraic Riccati equation, perturbation
analysis, condition number, statistical condition estimation.\\

}


\section{Introduction}\label{secintro}


Algebraic Riccati equations arise in optimal control problems in continuous-time or discrete-time. The theory, applications, and numerical methods for solving the equations can be found in
\cite{AM79,BLW91,LR95,PLD94,SW77} and references therein. The continuous-time algebraic Riccati equation (CARE) is given in the form:
\begin{equation}\label{careorigin}
Q+A^H X+XA-XBR^{-1}B^H X=0,
\end{equation}
where $X$ is the unknown matrix, $A \in \C^{n\times n}$, $B\in \C^{n\times m}$, $A^H$ denotes the conjugate transpose of $A$,  and $Q,\, R$ are $n\times n$ Hermitian matrices with $Q$ being positive semi-definite (p.s.d.) and $R$ being positive definite. The discrete-time algebraic Riccati equation (DARE) is given in the form:
\begin{equation}\label{dareorigin}
Y-A^H YA+A^H YB(R+B^H YB)^{-1}B^H YA-C^H C=0,
\end{equation}
where $Y$ is the unknown matrix,
$A \in \C^{n\times n}$, $B\in \C^{n\times m}$,
$C\in \C^{r\times n}$, and $R\in \C^{m\times m}$ with
$R$ being Hermitian positive definite.

For the complex CARE (\ref{careorigin}), let
$G=BR^{-1}B^H$, then it has the
simplified form
\begin{equation}\label{care}
Q+A^HX+XA-XGX=0,
\end{equation}
where $Q,G$ are Hermitian and p.s.d.
For the complex DARE (\ref{dareorigin}),
let $Q=C^HC$ and $G=BR^{-1}B^H$, then it has the
simplified form
\begin{equation}\label{dare}
Y-A^HY(I+GY)^{-1}A-Q=0,
\end{equation}
where $Q,G$ are Hermitian and p.s.d.
In particular, when $A$, $Q$ and $G$ are real matrices,
the real CARE becomes
\begin{equation}\label{Rcare}
Q+A^TX+XA-XGX=0,
\end{equation}
and the real DARE has the form
\begin{equation}\label{Rdare}
Y-A^TY(I+GY)^{-1}A-Q=0.
\end{equation}

The existence and uniqueness of the solution is essential
for perturbation analysis. Before making
appropriate assumptions on the coefficient matrices
necessary for the existence and uniqueness of Hermitian
and p.s.d. stabilizing solution, we need
some notions of stability, which play an important role in
the study of the algebraic Riccati equations. An $n\times n$ matrix
$M$ is said to be c-stable if all of its eigenvalues lie in the
open left-half complex plane, and $M$ is said to be d-stable if
its spectral radius $\rho(M)<1$.
Then to ensure the existence and uniqueness of the solution,
we assume that $(A,G)$ in the CARE
(\ref{care}) is a c-stabilizable pair, that is, there is a matrix
$K \in \C^{n\times n}$ such that the matrix $A-GK$ is c-stable,
and that $(A,Q)$ is a c-detectable pair, that is, $(A^T,Q^T)$ is
c-stabilizable. It is known \cite{Byers85,Laub79} that under
these conditions there exists a unique Hermitian and p.s.d.
solution $X$ for the CARE (\ref{care}) and the matrix $A-GX$ is
c-stable. Similarly, for the DARE, we assume that
$(A,B)$ in the DARE (\ref{dareorigin}) is a d-stabilizable pair,
that is, if $\omega^TB=0$ and $\omega^TA=\lambda \omega^T$ hold
for some constant $\lambda$, then $|\lambda|<1$ or $\omega=0$, and
that $(A,C)$ is a d-detectable pair, that is, $(A^T,C^T)$ is
d-stabilizable. It is known \cite{AM79,GKL92,KPC93} that under
these conditions there exists a unique Hermitian and p.s.d.
solution $Y$ for the DARE (\ref{dare}), and the matrix
$(I+GY)^{-1}A$ is d-stable, i.e., all the eigenvalues of
$(I+GY)^{-1}A$ lie in the open unit disk.

Matrix perturbation analysis concerns the sensitivity of the
solution to the perturbations in the data of a problem. A condition
number is a measurement of the sensitivity. Liu studied mixed and
componentwise condition numbers of nonsymmetric
algebraic Riccati equation in \cite{Liu}. For the
perturbation analysis of the CARE (\ref{care}) or DARE
(\ref{dare}), we refer papers \cite{Byers85,KH90,KLW90,GKL92,KPC93} and
their references therein. Sun~\cite{Sun02} defined
the structured normwise condition numbers
for CARE and DARE and showed that the expressions of structured
normwise condition numbers are the same as their
unstructured counterparts for both real
and complex cases. Later, Zhou et al. \cite{zhou} performed componentwise
perturbation analyses of CARE and DARE and obtained the exact expressions
for mixed and componentwise condition numbers defined in \cite{GK93}
for the real case. However, in their paper, the perturbations on
$Q$ and $G$ are general (unstructured).
In this paper,
we perform a structured perturbation analysis, define the
structured normwise, mixed and componentwise condition numbers
for complex CARE and DARE, and derive their expressions using
the Kronecker product \cite{Graham81}.
Specifically, we assume that the perturbation $\Delta G$
($\Delta Q$) has the same structure as $G$ ($Q$).
Furthermore, in the complex
case, we separate the real part and the imaginary part. Thus the real part
of $G$ ($Q$) or $\Delta G$ ($\Delta Q$) is symmetric and the
imaginary part of $G$ ($Q$) or $\Delta G$ ($\Delta Q$) is
skew-symmetric. In our analysis, we exploit the structure and
consider the perturbations on the real part and the imaginary part
separately. In contrast, the analysis in \cite{Sun02} considers
the perturbation on a complex matrix as whole. Apparently,
separating real and imaginary parts gives more precise results.

Efficiently estimating the condition
of a problem is one of the most fundamental topics in numerical analysis.
Together with the knowledge of backward error, a good condition
estimate can provide an estimate for the accuracy of the computed solution.
Although the expressions of the condition numbers derived in
\cite{Sun02}, \cite{zhou} and this paper are explicit, they
involve the solution matrix and require extensive computation,
especially for large size problems.
As pointed out in \cite[Page 260]{Sun02}, practical
algorithms for estimating the condition numbers of
the algebraic Riccati equations are worth studying.
In this paper, we present a statistical method for practically
estimating the structured normwise, mixed and componentwise condition numbers
for CARE and DARE by applying the small sample condition estimation
method (SCE) \cite{KenneyLaub_SISC94}. 

The SCE, proposed by Kenny and Laub~\cite{KenneyLaub_SISC94},
is an efficient method for estimating the condition numbers for
linear systems \cite{KenneyLinear,laublinear}, linear least squares problems \cite{KenneyLS}, the Tikhonov regularization problem \cite{DWQ}, the total least squares problem \cite{DWX}, eigenvalue problems \cite{laubeig},
roots of polynomials \cite{laubroot}, etc. Diao et al. \cite{DSW,DXW,DYC} applied
the SCE to the (generalized) Sylvester equations.
Wang et al. \cite {Wang16} considered the mixed and componentwise
condition numbers for the spectral projections, generalized spectral
projections and sign functions for matrices and regular matrix pairs
and derived explicit expressions of the condition numbers, which
improved some known results of the normwise type and revealed the
structured perturbations. Also, they applied the SCE to these problems
to efficiently estimate the condition numbers.
Wang et al. \cite{Wang15} studied the normwise, mixed
and componentwise condition numbers for the following general
nonlinear matrix equation
$
X+A^H F(X) A=Q,
$
where $A$ is an $n$-by-$n$ square matrix, $Q$ an $n$-by-$n$
positive definite matrix, $X$ the unknown $n$-by-$n$ positive
semi-definite matrix, and $F$ a differentiable mapping from the set
of $n$-by-$n$ positive semi-definite matrices to the set of
$n$-by-$n$ matrices. They derived corresponding explicit condition numbers and gave
their statistical estimations with high reliability based on the
SCE and a probabilistic spectral norm estimator. Differing from
their algorithms, our methods produce estimated condition matrices
instead of single condition numbers, that is, the entries of the
condition matrices produced by our algorithms are the structured
normwise or componentwise condition numbers of the corresponding
entries of the solution matrices. Thus these condition matrices
are more informative and precise about the conditioning of the solution.

Throughout this paper we adopt the following notations:

\begin{itemize}
\item
$\C^{m\times n}$ ($\R^{m\times n}$)  denotes the set of
complex (real) $m \times n$ matrices;
$\HM^{n\times n}$ the set of $n \times n$ Hermitian matrices;
$\SM^{n\times n}$ the set of $n \times n$ symmetric matrices;
$\SK^{n\times n}$ the set of $n \times n$ skew-symmetric matrices.
\item
$A^T$ denotes the transpose of $A$;
$A^H$ the complex conjugate and transpose of $A$;
$A^{\dagger}$ the Moore-Penrose inverse of $A$;
$I$ the identity matrix;
$0$ the zero matrix;
$\Re(A)$ ($\Im (A)$) is the real (imaginary)
part of a complex matrix $A$. The matrix $\diag(A,B)$ denotes a
block diagonal matrix with $A$ and $B$ being its diagonal.
\item
$e_i$ denotes the $i$th column of $I$.
\item
The mapping $\sym(\cdot)$:
$\SM^{n\times n} \rightarrow \R^{n(n+1)/2}$ maps a symmetric matrix
$A=[a_{ij}] \in \SM^{n\times n}$ to a $(n(n+1)/2)$-vector:
$$
[a_{11}, ..., a_{1n},a_{22}, ..., a_{2n}, ...,a_{n-1,n-1},
a_{n-1,n},a_{nn}]^T.
$$
\item
The mapping $\ske(\cdot)$:
$\SK^{n\times n} \rightarrow \R^{n(n-1)/2}$ maps a skew-symmetric
matrix $A=[a_{ij}] \in \SK^{n\times n}$ to the $(n(n-1)/2)$-vector:
$$
[a_{12},..., a_{1n},a_{23}, ..., a_{2n}, ...,a_{n-2,n-1},
a_{n-2,n} a_{n-1,n}]^T.
$$
\item
$A \succ 0$ ($A \succeq 0$) means that $A$ is positive
definite (positive semi-definite).
\item
$\|~\|_F, \|~\|_2$ and $\|~\|_\infty$ are
the Frobenius norm, the spectral norm and infinity norm respectively.
For $A \in \C^{m\times n}$, $\|A\|_{\max}=\max_{ij}|a_{ij}|$.
\item
$A \otimes B = [a_{ij}B]$ is the Kronecker product
of $A = [a_{ij}]$ and matrix $B$ and
$\vect(A)$ is the vector defined by
$\vect(A)=[a_1^T, ...,a_n^T]^T \in \C^{mn}$;
$\Pi$ is an $n^2\times n^2$
permutation matrix, such that, for an $n\times n$ real matrix $A$,
$\vect(A^T) = \Pi \vect(A)$. For more properties of
the Kronecker product and $\vect$
operation, see \cite{Graham81}.
\item
$|A|\leq |B|$ means $|a_{ij}|\leq |b_{ij}|$ for $A,\, B\in \C^{m\times n}$;
$A \oslash B$ is the componentwise division of matrices $A$ and $B$ of
the same dimensions. In our context, it is used for componentwise
relative error. So, when $b_{ij}= 0$, we assume its absolute
error $a_{ij} = 0$ and set $(A \oslash B)_{ij}=0$.
\end{itemize}

The rest of the paper is organized as follows.
In Section~\ref{secnorm}, we present our structured perturbation
analyses and expressions of the structured normwise
condition numbers of the CARE (\ref{care}) and the DARE (\ref{dare}).
The expressions of the structured
mixed and componentwise condition numbers are derived in
Section~\ref{seccomponentnorm}. In Section~\ref{secsce},
by applying the small sample condition estimation method,
we propose our structured sensitivity estimation methods for the
problems of solving the CARE and DARE. Our numerical experiment
results are demonstrated in Section~\ref{secnumex}.
Finally, Section~\ref{secconclude} concludes this paper.

\section{\bf Structured Normwise Condition Numbers}\label{secnorm}

In this section, using the Kronecker product, we first present a
structured perturbation analysis of the CARE (\ref{care}) and derive
expressions of the corresponding structured
normwise condition number. In a similar way, a
structured perturbation analysis of the DARE (\ref{dare}) can
be performed and
the corresponding structured normwise condition number
can be obtained.

\subsection{CARE}

Let $\Delta A \in \C^{n\times n}$, $\Delta Q \in \HM^{n\times n}$,
and $\Delta G\in \HM^{n\times n}$ be the perturbations to the
data $A$, $Q$, and $G$ respectively. Notice that the perturbations
$\Delta Q$ and $\Delta G$ are also Hermitian.
From Theorem~3.1 and its proof
in \cite{Sun98SIAM}, for $Q,\, G \succeq 0$ and sufficiently
small $\|[\Delta A, \Delta Q, \Delta G]\|_F$, there is a unique
Hermitian p.s.d. matrix $\tilde X$ such
that  $\tilde A-\tilde G \tilde X$ is c-stable, where
$\tilde A=A+\Delta A$, $\tilde G=G+\Delta G$, and
\begin{equation}\label{perturbequ}
\tilde X\tilde G\tilde X- \tilde X \tilde A-\tilde A^H\tilde X-\tilde Q =0,
\end{equation}
where $\tilde Q=Q+\Delta Q$.
To perform a perturbation analysis, we
define a linear operator ${\bf L}:\
\HM^{n\times n} \rightarrow \HM^{n\times n}$ by
\begin{equation} \label{eqnLW}
{\bf L} W=(A-GX)^H W+W(A-GX), \quad W\in \HM^{n\times n}.
\end{equation}
Since the matrix $A-GX$ is c-stable, the operator $\bf L$ is
invertible. Denote $\Delta X=\tilde X-X$, then
\begin{equation}\label{eq:deltaXCARE}
\Delta X = -{\bf L}^{-1} (\Delta Q + X\Delta A +
\Delta A^H X - X\Delta G X)  + O\left(\|[\Delta A, \Delta Q, \Delta G]\|_F^2\right),
\end{equation}
as $\|[\Delta A, \Delta Q, \Delta G]\|_F \rightarrow 0$.
In the first order approximation of (\ref{eq:deltaXCARE}),
$\Delta X$ is the solution
to the continuous Lyapunov equation:
\begin{equation}
 (A - GX)^H \Delta X + \Delta X (A - GX)
 =  -\Delta Q - X\Delta A - \Delta A^H X + X\Delta G X .
\label{eq:Lyapunov}
\end{equation}
Numerical methods for solving the Lyapunov equations
can be found in \cite{BIM12}.

For the complex CARE (\ref{care}), in addition to exploiting the
symmetry structure of $G$ and $Q$, to make our analysis precise,
we separate the real part and the imaginary part. Since the real
part of a Hermitian matrix is symmetric and the imaginary part
is skew-symmetric, we introduce the function
\begin{align*}
&\varphi :\ \R^{4n^2}\rightarrow \C^{n^2}\\
&[\vect(\Re(A))^T,\vect(\Im(A))^T,\sym(\Re(G))^T,\ske(\Im(G))^T,\\
&\  \sym(\Re(Q))^T, \ske(\Im(Q))^T ]^T \mapsto \vect(X) ,
\end{align*}
which maps the structured data vector to the solution vector.
As we can see, our data vector exploits the structure and
separates the real and imaginary parts.
Applying the condition number theory of Rice \cite{Rice}
to the above mapping,
we define the structured normwise condition numbers:
\begin{equation}\label{defnm}
\kappa_i (\varphi) = \lim_{\epsilon \rightarrow 0}
\sup_{\substack{\eta_i \le\epsilon \\
\Delta A\in \C^{n\times n},\, \Delta Q\in \HM^{n\times n},\,
\Delta G\in \HM^{n\times n}} \atop {Q+\Delta Q \succeq 0,\,
G+\Delta G \succeq 0    }}  \frac{\|\Delta X\|_F}
{\epsilon \|X\|_F},
\end{equation}
where
\begin{eqnarray}
\eta_1&=&\Big\|\Big[ \frac{\|\Re(\Delta A)\|_F}{\delta_1},\,
\frac{\|\Im(\Delta A)\|_F}{\delta_2},
\frac{\|\sym(\Re(\Delta G))\|_2}{\delta_3},\,
\frac{\|\ske(\Im(\Delta G))\|_2}{\delta_4},
\nonumber \\
&&\frac{\|\sym(\Re(\Delta Q))\|_2}{\delta_5},\,
\frac{\|\ske(\Im(\Delta Q))\|_2}{\delta_6} \Big]\Big\|_2, \nonumber \\
\eta_2&=&\max \Big\{ \frac{\|\Re(\Delta A)\|_F}{\delta_1},\,
\frac{\|\Im(\Delta A)\|_F}{\delta_2}, \frac{\|\sym(\Re(\Delta G))\|_2}{\delta_3},\,
\frac{\|\ske(\Im(\Delta G))\|_2}{\delta_4}, \nonumber \\
&&\frac{\|\sym(\Re(\Delta Q))\|_2}{\delta_5},\,
\frac{\|\ske(\Im(\Delta Q))\|_2}{\delta_6} \Big\},\label{defDelta}
\end{eqnarray}
and the parameters $\delta_i > 0 \,(i=1,...,6)$ are given.
Generally, they are respectively chosen to be the functions of
$\|\Re(A)\|_F$, $\|\Im(A)\|_F$, $\|\sym(\Re(G))\|_2$,
$\|\ske(\Im(G))\|_2$, $\|\sym(\Re(Q))\|_2$, and $\|\ske(\Im(Q))\|_2$.
Here, we set $\delta_1=\|\Re(A)\|_F$, $\delta_2=\|\Im(A)\|_F$,
$\delta_3=\|\sym(\Re(G))\|_2$, $\delta_4=\|\ske(\Im(G))\|_2$,
$\delta_5=\|\sym(\Re(Q))\|_2$ and $\delta_6=\|\ske(\Im(Q))\|_2$.

Now, we derive an explicit expression of $\kappa_1 ( \varphi )$
and an upper bound for $\kappa_2 ( \varphi )$ in (\ref{defnm}).
First, we present a matrix-tensor representation of the operator
$\mathbf{L}$. Applying the identity
\begin{equation} \label{eqnIdentity}
\vect (UVW) = (W^T \otimes U) \vect (V),
\end{equation}
to the vectorized (\ref{eqnLW})
\[
\vect ({\bf L} W) =
\vect ((A-GX)^H W+W(A-GX) ),
\]
we get the matrix
\begin{equation} \label{eqnZ}
Z=I_n\otimes (A-GX)^H+(A-GX)^T \otimes I_n ,
\end{equation}
which transforms $\vect (W)$ into $\vect ( \mathbf{L} W)$,
as a matrix representation of the linear operator $\mathbf{L}$.
Since $\bf L$ is invertible, $Z$ is also invertible.

Since our structured data vector exploits the symmetry structure,
to convert it back to original vector, we introduce the matrices
$\mathcal{S}_1$ and $\mathcal{S}_2$ as follows.
For an $n \times n$ symmetric matrix $J$, $\mathcal{S}_1$
is the $n^2 \times n(n+1)/2$ matrix such that
\[
\vect (J) = \mathcal{S}_1 \sym (J) .
\]
That is, $\mathcal{S}_1$ expands the $n(n+1)/2$-vector
$\sym (J)$ to the $n^2$-vector $\vect (J)$ by copying
its elements. The $n^2 \times n(n-1)/2$ matrix $\mathcal{S}_2$
is defined by
\[
\vect (K) = \mathcal{S}_2 \ske (K) ,
\]
where $K$ is an $n \times n$  skew-symmetric matrix.

Then, applying (\ref{eqnIdentity}) to
(\ref{eq:deltaXCARE}), we get
\begin{align}
 &\vect(\Delta X)
= \mbox{} -Z^{-1} [-(X^T\otimes X)\vect(\Delta G) +
(I_n\otimes X) \vect(\Delta A) \nonumber \\
&  \mbox{} +(X^T\otimes I_n)\vect(\Delta A^H) +
\vect(\Delta Q) ]+O(\|[\Delta A, \Delta Q, \Delta G]\|_F^2) \nonumber \\
&= \mbox{} -Z^{-1} [-(X^T\otimes X)(\cals_1\sym(\Re(\Delta G))  + \bfi \, \cals_2 \ske(\Im(\Delta G))) \nonumber \\
&  \mbox{}+
(I_n\otimes X)(\vect(\Re(\Delta A)) +\bfi \, \vect(\Im(\Delta A)))
+(X^T\otimes I_n) \Pi (\vect(\Re(\Delta A)) \nonumber \\
&  \mbox{} -\bfi \, \vect(\Im(\Delta A))) +
\cals_1 \sym(\Re(\Delta Q))  +\bfi \, \cals_2\ske(\Im(\Delta Q))] +
O(\|[\Delta A, \Delta Q, \Delta G]\|_F^2)\nonumber \\
&\to  -Z^{-1}[(I_n\otimes X) + (X^T\otimes I_n)\Pi, \bfi \, ((I_n\otimes X) - (X^T\otimes I_n)\Pi), \nonumber\\
&  \ -(X^T\otimes X)\cals_1,\
-\bfi \, (X^T\otimes X)\cals_2,\ \cals_1, \ \bfi \, \cals_2]
\cdot \Delta \label{Zvec} ,
\end{align}
as $\|[\Delta A, \Delta Q, \Delta G]\|_F \to 0$,
where $\bfi=\sqrt{-1}$ and
\begin{eqnarray}
\Delta &=&
[\vect(\Re(\Delta A))^T,\ \vect(\Im(\Delta A))^T,\ \sym(\Re(\Delta G))^T,\ \ske(\Im(\Delta G))^T,\nonumber \\
& & \ \sym(\Re(\Delta Q))^T,\ \ske(\Im(\Delta Q))^T]^T
\label{eq:Delta all}
\end{eqnarray}
is the structured data perturbation vector.

Denoting
\begin{eqnarray*}
M_A &=& [(I_n\otimes X) + (X^T\otimes I_n)\Pi,\ \bfi((I_n\otimes X) -(X^T\otimes I_n)\Pi)]
\end{eqnarray*}
corresponding to
$[\vect(\Re(\Delta A))^T,\ \vect(\Im(\Delta A))^T]^T$, 
$$
M_G = [-(X^T\otimes X)\cals_1,\  -\bfi(X^T\otimes X)\cals_2]
$$
corresponding to
$[\sym(\Re(\Delta G))^T,\ \ske(\Im(\Delta G))^T]^T$,
$$
M_Q = [\cals_1, \ \bfi \, \cals_2]
$$
corresponding to
$[\sym(\Re(\Delta Q))^T,\ \ske(\Im(\Delta Q))^T]^T$, and
$\calm = [M_A \ M_G \ M_Q]$,
and using the definition of the directional derivative,
we have the following lemma.

\begin{lemma} \label{lemma:care}
Using the above notations, the directional derivative
$\cald \varphi(X)$ of
$\varphi$ 
with respect to $\Delta$ (\ref{eq:Delta all}) is given by
\begin{equation}\label{eq:derative varphi}
\cald \varphi(X)\Delta=-Z^{-1} \calm \cdot \Delta.
\end{equation}
\end{lemma}

Finally, the following theorem gives an explicit expression
of $\kappa_1(\varphi )$ and an
upper bound for $\kappa_2(\varphi )$.

\begin{theorem}\label{normthm}
Using the notations given above, the expression and upper bound
for the normwise number of the complex CARE
(\ref{care}) are
\begin{eqnarray}
\kappa_1(\varphi)&= & \frac{\|Z^{-1}\calm D\|_2}{\|X\|_F} \label{expnm-1},\\
\kappa_2(\varphi)&\leq & \kappa_U(\varphi):=
\min\left\{\sqrt{6}\kappa_1(\varphi), \alpha_c/\|X\|_F\right\},\label{expnm-2}
\end{eqnarray}
where
\begin{align}
D &= \diag \left([\delta_1 I_{n^2},\ \delta_2 I_{n^2},\
\delta_3 I_{n(n+1)/2},\delta_4I_{n(n-1)/2},\ \delta_5 I_{n(n+1)/2},\
\delta_6I_{n(n-1)/2}]\right.) \label{eq:D}
\end{align}
and
\begin{align*}
\alpha_c &= \delta_1 \|Z^{-1}[(I_n\otimes X) + (X^T\otimes I_n)\Pi]\|_2 + \delta_2\|Z^{-1}[(I_n\otimes X)-(X^T\otimes I_n)\Pi ]\|_2 \\
& \mbox{}+ \delta_3 \|Z^{-1} (X^T\otimes X)\cals_1\|_2+ \delta_4 \|Z^{-1} (X^T\otimes X)\cals_2 \|_2 + \delta_5 \|Z^{-1}\cals_1\|_2+\delta_6\| Z^{-1}\cals_2\|_2.
\end{align*}
\end{theorem}
{\bf Proof.}
Introducing the positive parameters $\delta_i$, $i=1,...,6$,
into (\ref{eq:derative varphi}), we get
\[
\cald\varphi(X)\Delta= -Z^{-1}\calm D D^{-1}\Delta.
\]
From the definition \eqref{defnm}, we know that
\begin{align}\label{deltaxbound-1}
\kappa_1(\varphi)&=\max_{\eta_1 \leq \epsilon}
\frac{\|-Z^{-1}\calm D D^{-1}\Delta\|_2}{\|X\|_F}=\max_{\|D^{-1}\Delta \|_2 \leq 1}
\frac{\|-Z^{-1}\calm D D^{-1}\Delta \|_2}{\|X\|_F}\nonumber \\
&=\frac{\|Z^{-1}\calm D  \|_2}{\|X\|_F}.
\end{align}
The last equality holds because $\Delta$ can vary freely.

Because $\|D^{-1}\Delta \|_2\leq \sqrt{6}\eta_2$, it is easy to see that
\begin{eqnarray*}
\|\Delta X\|_F &\approx& \|Z^{-1}\calm D D^{-1}\Delta\|_2 \le\|Z^{-1}\calm D \|_2\|D^{-1}\Delta\|_2\leq \sqrt 6 \|Z^{-1}\calm D \|_2\eta_2,
\end{eqnarray*}
which proves $\kappa_2(\varphi)\leq  \sqrt{6}\kappa_1(\varphi)$.
On the other hand, since
\begin{eqnarray*}
& & \vect(\Delta X) = -\delta_1 Z^{-1}[(I_n\otimes X) +(X^T\otimes I_n)\Pi]
\frac{\vect(\Re(\Delta A)))}{\delta_1} \\
& & \mbox{} -\bfi \delta_2Z^{-1}[(I_n\otimes X) -(X^T\otimes I_n)\Pi ]
\frac{\vect(\Im(\Delta A)))}{\delta_2} \\
& & \mbox{} +\delta_3 Z^{-1} (X^T\otimes X)\cals_1
\frac{\sym(\Re(\Delta G))}{\delta_3} \mbox{} +\bfi\delta_4 Z^{-1} (X^T\otimes X)\cals_2
\frac{\ske(\Im(\Delta G))}{\delta_4} \\
& & \mbox{} -\delta_5 Z^{-1}\cals_1
\frac{\sym(\Re(\Delta Q))}{\delta_5}\mbox{} -\bfi\delta_6 Z^{-1}\cals_2
\frac{\ske(\Im(\Delta Q))}{\delta_6},
\end{eqnarray*}
it is easy to see that $\|\Delta X\|_F \leq \alpha_c \eta_2$. \qed

In \cite{Sun02}, Sun defined the structured normwise
condition number for the complex CARE:
\begin{equation}\label{defSunCARE}
\kappa_{\rm Sun}^{\rm CARE} = \lim_{\epsilon \to 0}
\sup_{\substack{\theta \le\epsilon \\
\Delta A \in \C^{n\times n},\,
\Delta G, \, \Delta Q \in \HM^{n\times n}} \atop
{G + \Delta G, \, Q + \Delta Q \succeq 0}}
\frac{\|\Delta X\|_F}{\epsilon \|X\|_F},
\end{equation}
where $ \theta=\left\|\left[\Delta A / \mu_1,\
\Delta G / \mu_2, \ \Delta Q / \mu_3 \right]\right\|_F$,
and the parameters $\mu_i$, $i=1,2,3$, are positive.

Differently from the above Sun's definition,
our definition separates the relative perturbations in the
real and imaginary parts of a data matrix. Apparently,
it is more precise than the Sun's definition. Moreover,
it is more realistic, since, in computation, the real and
imaginary parts of a complex matrix are stored and computed
separately.

From the definitions and expressions of \eqref{defnm} and \eqref{defSunCARE},
if we choose $\delta_1=\delta_2=\mu_1$, $\delta_3=\delta_4=\mu_2$,
and $\delta_5=\delta_6=\mu_3$, we can show that
$$
\kappa_1(\varphi) \leq \max\{\|\cals_1\|_2,\, \|\cals_2\|_2\}\cdot\kappa_{\rm Sun}^{\rm CARE}.
$$
Note that each row of $\cals_1$ or $\cals_2$ is $e_i^T$ for
some $i$, that is, a row of the identity matrix. 

Similarly to the complex case,
for the real CARE (\ref{Rcare}), we introduce the mapping
\begin{eqnarray*}
\varphi^\Re &:& \R^{n^2+n(n+1)}\rightarrow \R^{n^2} \\
& & (\vect(A)^T,\sym(G))^T,\sym(Q)^T)^T\mapsto \vect(X)
\end{eqnarray*}
and define the structured normwise condition
number for the real CARE \eqref{Rcare}:
\[
\kappa (\varphi^\Re) = \lim_{\epsilon \to 0}
\sup_{\substack{\delta \le \epsilon \\
\Delta A \in \R^{n\times n},\,
\Delta G, \, \Delta Q \in \SM^{n\times n}}\atop
{G + \Delta G, \, Q + \Delta Q \succeq 0}}
\frac{\|\Delta X\|_F}{\epsilon \|X\|_F},
\]
where
\begin{equation} \label{defDeltaR}
\delta = \max \left\{
\frac{\|\Delta A\|_F}{\delta_1},\,
\frac{\|\sym(\Delta G)\|_2}{\delta_2},\,
\frac{\|\sym(\Delta Q)\|_2}{\delta_3} \right\},
\end{equation}
with the positive parameters $\delta_i$, $i=1,2,3$.
We usually choose $\delta_1=\|A\|_F$, $\delta_2=\|\sym(G)\|_2$ and
$\delta_3=\|\sym(Q)\|_2$.

Using the above notations, we can derive the following
upper bound for the normwise condition number of the real CARE
(\ref{Rcare}):
\[
\kappa(\varphi^\Re) \leq \kappa_U(\varphi^\Re):=
\min \left\{\sqrt 3 \frac{\|Z_1^{-1}\calm_1 D_1\|_2}
{\|X\|_F}, \, \frac{\beta}{\|X\|_F} \right\},
\]
where
\begin{align*}
Z_1&=I_n\otimes (A-GX)^T+(A-GX)^T \otimes I_n,\\
\calm_1&=\Big[I_n\otimes X+(X^T\otimes I_n)\Pi,
~-(X^T\otimes X)\cals_1,~\cals_1\Big],\\
D_1&=\diag\left([\delta_1 I_{n^2},\ \delta_2 I_{n(n+1)/2},\
\delta_3 I_{n(n+1)/2}]\right.),
\end{align*}
and
\begin{eqnarray*}
\beta &=& \delta_1 \|Z_1^{-1}[(I_n\otimes X) + (X\otimes I_n)\Pi]\|_2\mbox{} +\delta_2\|Z_1^{-1}(X\otimes X)\cals_1\|_2+
\delta_3\|Z_1^{-1}\cals_1\|_2.
\end{eqnarray*}

In \cite{zhou}, Zhou et al. defined the following
unstructured normwise condition number for the real CARE:
$$
\kappa_1 (\varphi^\Re) = \lim_{\epsilon \rightarrow 0}
\sup_{\delta \le\epsilon} \frac{\|\Delta X\|_F}{\epsilon \|X\|_F},
$$
where
\begin{eqnarray*}
\delta&=&\max \left\{ \frac{\|\Delta A\|_F}{\delta_1},\,
\frac{\|\Delta Q\|_F}{\delta_2},\,
\frac{\|\Delta G\|_F}{\delta_3}\right\},
\end{eqnarray*}
and derived the upper bound
$$
\kappa_1(\varphi^\Re) \leq \kappa_1^U (\varphi^\Re) :=
\min\left\{\sqrt{3}\,\frac{\|Z_1^{-1}S_1\|_2}{\|X\|_F}, \
\frac{\beta_c}{\| X \|_F} \right\},
$$
where
\begin{eqnarray*}
S_1 &=& [-(I_n\otimes X)-(X\otimes I_n) \Pi,  X\otimes X, \ -I_{n^2}] \,
\diag([\delta_1 I_{n^2},\ \delta_2 I_{n^2},\ \delta_3I_{n^2}])
\end{eqnarray*}
and
\begin{eqnarray*}
\beta_c &=& \delta_1 \|Z_1^{-1}[I_n\otimes X+(X\otimes I_n)\Pi]\|_2  \mbox{} +\delta_2\|Z_1^{-1}(X\otimes X)\|_2+\delta_3\|Z_1^{-1}\|_2.
\end{eqnarray*}

By setting $\delta_1=\|A\|_F$, $\delta_2=\|Q\|_F$ and
$\delta_3=\|G\|_F$, we can prove that
$$
\kappa_U(\varphi^\Re) \leq \|\cals_1\|_2\kappa_1^U(\varphi^\Re).
$$
However,
our numerical experiments show that the difference between
our $\kappa_U(\varphi^\Re)$ and  $\kappa_1^U(\varphi^\Re)$
is marginal.

\subsection{DARE}

Following the structured perturbation analysis of CARE,
for the complex DARE (\ref{dare}), we define
a linear operator
${\bf L}:\HM^{n\times n} \rightarrow \HM^{n\times n}$ by
$$
{\bf L} M = M-[(I_n+G Y)^{-1}A]^H M (I_n+G Y)^{-1}A,
$$
for $M\in \HM^{n\times n}$.
Since the matrix $(I_n+G Y)^{-1}A$ is d-stable,
the operator $\bf L$ is invertible.
Let
\begin{equation} \label{perturbequdare}
\tilde Y-\tilde A^H \tilde Y(I_n+\tilde G \tilde  Y )^{-1}\tilde A-
\tilde Q = 0,
\end{equation}
be the perturbed DARE, where $\tilde A=A+\Delta A$,
$\tilde G=G+\Delta G$, and $\tilde Q=Q+\Delta Q$, then
for sufficiently small $\|[\Delta A, \Delta Q, \Delta G]\|_F$,
there is a unique Hermitian and p.s.d. solution
$\tilde Y$ for the perturbed equation (\ref{perturbequdare})
and the change $\Delta Y=\tilde Y-Y$ in the solution
is given by
\begin{equation}
\Delta Y \to {\bf L}^{-1}[\Delta Q+(A^HYW)\Delta A+
\Delta A^H (YWA)  \mbox{} -(A^HYW)\Delta G (YWA)], \label{eq:delta XDARE}
\end{equation}
as $\|[\Delta A, \Delta Q, \Delta G]\|_F \to 0$.

Denote $W=(I_n+G Y)^{-1}$. In this case, the
matrix-tensor representation of the linear operator $\bf L$ is
\begin{equation} \label{eqnT}
T=I_n-(A^TW^T)\otimes (A^H W^H),
\end{equation}
which is invertible, since $\bf L$ is invertible.

For the structured perturbation analysis of the
complex DARE (\ref{dare}), exploiting the symmetry
structure and separating the real and imaginary parts,
we define the mapping
\begin{eqnarray*}
\psi &:& \R^{4n^2} \rightarrow \C^{n^2}\\
& & [\vect(\Re(A))^T,\ \vect(\Im(A))^T,\ \sym(\Re(G))^T, \\
& & \quad \ske(\Im(G))^T, \ \sym(\Re(Q))^T, \ \ske(\Im(Q))^T ]^T \\
& & \mapsto \vect(Y).
\end{eqnarray*}
Similarly to
(\ref{Zvec}), dropping the second and higher order terms, we have
\begin{eqnarray*}
& & T\vect(\Delta Y) \approx [(I_n\otimes(A^H YW)) +((A^TW^TY^T)\otimes I_n)\Pi,\\
& & \ \bfi ((I_n\otimes(A^H YW)) -((A^TW^TY^T)\otimes I_n)\Pi ), \ -(A^TW^TY^T)\otimes(A^HYW)\cals_1,\\
& & \ -\bfi ((A^TW^TY^T)\otimes(A^HYW)) \cals_2,
\ \cals_1,\ \bfi \cals_2] \cdot \Delta,
\end{eqnarray*}
where the data perturbation vector $\Delta$ is defined in
\eqref{eq:Delta all}.
Denoting
\begin{align*}
N_A &=[(I_n\otimes(A^H YW)) +((A^TW^TY^T)\otimes I_n)\Pi, \bfi ((I_n\otimes(A^H YW)) -((A^TW^TY^T)\otimes I_n)\Pi )], \\
N_G &=[-(A^TW^TY^T)\otimes(A^HYW)\cals_1, -\bfi ((A^TW^TY^T)\otimes(A^HYW)) \cals_2], \\
N_Q &=[\cals_1,\ \bfi \cals_2] ,
\end{align*}
and $\caln = [N_A \ N_G \ N_Q]$, we have the following lemma.

\begin{lemma}\label{lemma:dare}
With the above notations, the directional derivative
$\cald\psi (Y)$ of
$\psi$ 
with respect to
$\Delta$ is
$$
\cald\psi (Y)\Delta_1=T^{-1}\caln\cdot \Delta.
$$
\end{lemma}

We then define the structured normwise
condition numbers for solving the complex DARE \eqref{dare}:
\[
\kappa_i (\psi)=
\lim_{\epsilon \rightarrow
0}\sup_{\substack{\eta_i\le\epsilon \\ \Delta A\in \C^{n\times n},\, \Delta Q\in \HM^{n\times n},\, \Delta G\in \HM^{n\times n}}\atop {Q+\Delta Q \succeq 0,\, G+\Delta G \succeq 0    }}
\frac{\|\Delta
Y\|_F}{\epsilon \|Y\|_F},
\]
where $\eta_i$, $i=1,2$, are defined in (\ref{defDelta}).

Similarly to the proof of Theorem~\ref{normthm}, we can obtain an
explicit expression of $\kappa_1 ( \psi )$ and an upper bound for
$\kappa_2 (\psi )$ given in the following theorem.

\begin{theorem}\label{normthmdare}
Using the above notations, an explicit expression and
an upper bound for the structured normwise numbers of the complex DARE
(\ref{dare}) are
\begin{eqnarray*}
\kappa_1(\psi)&=& \frac{\|T^{-1}\caln D\|_2}{\|Y\|_F},\\
\kappa_2(\psi)&\leq& \kappa_U(\psi):=
\min\left\{\sqrt 6 \kappa_1(\psi),\alpha_d/\|Y\|_F\right\},
\end{eqnarray*}
where $D$ is defined in \eqref{eq:D} and
\begin{eqnarray*}
& & \alpha_d
= \delta_1\|T^{-1}((I_n\otimes(A^H YW)) +
((A^TW^TY^T)\otimes I_n)\Pi )\|_2\\
& & \mbox{} + \delta_2\left\|T^{-1}((I_n\otimes(A^H YW)) -
((A^TW^TY^T)\otimes I_n)\Pi)\right\|_2 \\
& & \mbox{} +\delta_3\|T^{-1}((A^TW^TY^T)\otimes(A^HYW)\cals_1 ) \|_2 \mbox{} +\delta_4\|T^{-1}((A^TW^TY^T)\otimes(A^HYW)\cals_2 ) \|_2  \\
& & \mbox{} +\delta_5\|T^{-1}\cals_1\|_2+\delta_6\|T^{-1}\cals_2\|_2.
\end{eqnarray*}
\end{theorem}

Using the parameter $\theta$ in \eqref{defSunCARE},
Sun \cite{Sun02} studied the structured normwise condition number
$\kappa_{\rm Sun}^{\rm DARE}$ for the complex DARE as follows
\begin{equation*}
\kappa_{\rm Sun}^{\rm DARE} = \lim_{\epsilon \to 0}
\sup_{\substack{\theta \le\epsilon \\
\Delta A \in \C^{n\times n},\,
\Delta G, \, \Delta Q \in \HM^{n\times n}} \atop
{G + \Delta G, \, Q + \Delta Q \succeq 0}}
\frac{\|\Delta Y\|_F}{\epsilon \|Y\|_F}.
\end{equation*}
Similarly to
the complex CARE case, we can prove that
$$
\kappa_1(\psi) \leq  \max\{\|\cals_1\|_2,\, \|\cals_2\|_2\}\cdot \kappa_{\rm Sun}^{\rm DARE},
$$
when we choose $\delta_1=\delta_2=\mu_1$,
$\delta_3=\delta_4=\mu_2$, and $\delta_5=\delta_6=\mu_3$.

For the real DARE (\ref{Rdare}), exploiting the
symmetry structure, we define the mapping
\begin{eqnarray*}
\psi^\Re &:& \R^{n^2+n(n+1)}\rightarrow \C^{n^2}\\
& & [\vect(A)^T,\sym(G)^T,\sym(Q)^T]^T \mapsto \vect(Y)
\end{eqnarray*}
and the structured normwise condition number:
\[
\kappa (\psi^\Re) = \lim_{\epsilon \rightarrow 0}
\sup_{\substack{\delta\le\epsilon \\
\Delta A\in \R^{n\times n},\,
\Delta G,\, \Delta Q\in \SM^{n\times n}} \atop
{G + \Delta G, \, Q + \Delta Q \succeq 0}}
\frac{\|\Delta Y\|_F}{\epsilon \|Y\|_F},
\]
where $\delta$ is defined in \eqref{defDeltaR}.

Using the above notations, we have the following
upper bound for the structured normwise condition
number of the real DARE:
\begin{equation} \label{expnm-3}
\kappa(\psi^\Re) \leq \kappa_U(\psi^\Re) :=\min\left\{
\frac{\sqrt{3}\, \|T_1^{-1}\caln_1 D_1\|_2}{\|Y\|_F},\,
\frac{\gamma}{\|Y\|_F} \right\},
\end{equation}
where
\begin{eqnarray*}
T_1 &=& I_n-(A^TW^T)\otimes (A^T W^T), \\
\caln_1 &=& [(I_n\otimes(A^T YW)) +((A^TW^TY^T)\otimes I_n)\Pi,\ -(A^TW^TY^T)\otimes(A^T YW)\cals_1,\cals_1 ],\\
D_1 &=& \diag ([\delta_1 I_{n^2},\ \delta_2 I_{n(n+1)/2},\
\delta_3 I_{n(n+1)/2}]),\\
\gamma&=&\delta_1\left\|T_1^{-1}[(I_n\otimes(A^T YW)) +
((A^TW^TY^T)\otimes I_n)\Pi]\right\|_2 \\
& & \mbox{} + \delta_2\left\|T_1^{-1} ((A^TW^TY^T)\otimes
(A^T YW)) \cals_1 \right\|_2  +\delta_3\left\|T_1^{-1} \cals_1 \right\|_2.
\end{eqnarray*}

As expected, our upper bound (\ref{expnm-3}) is an improvement of the
condition number in \cite[page 260]{Sun02}.

In \cite{zhou}, Zhou et al. defined the following
unstructured normwise condition number for the real DARE:
$$
\kappa_1 (\psi^\Re) = \lim_{\epsilon \rightarrow 0}
\sup_{\delta \le\epsilon} \frac{\|\Delta Y\|_F}{\epsilon \|Y\|_F},
$$
where
\begin{eqnarray*}
\delta&=&\max \left\{ \frac{\|\Delta A\|_F}{\delta_1},\,
\frac{\|\Delta Q\|_F}{\delta_2},\,
\frac{\|\Delta G\|_F}{\delta_3}\right\},
\end{eqnarray*}
and derived the upper bound
$$
\kappa_1(\psi^\Re) \leq \kappa_1^U (\psi^\Re) :=
\min\left\{\sqrt{3}\,\frac{\|T_1^{-1}P_1\|_2}{\|Y\|_F}, \
\frac{\gamma_d}{\| Y \|_F} \right\},
$$
where
\begin{eqnarray*}
P_1&=& \left[ (I_n\otimes(A^TYW)) + ((A^TW^TY)\otimes I_n) \Pi ,-((A^TW^TY)\otimes(A^TYW)),\, I_{n^2} \right] \\
& & \cdot \diag([\delta_1 I_{n^2},\ \delta_2 I_{n^2},\ \delta_3I_{n^2}])
\end{eqnarray*}
and
\begin{eqnarray*}
\gamma_d&=&\delta_1 \|T_1^{-1} ((I_n\otimes (A^TYW)) +
((A^TW^TY)\otimes I_n) \Pi )\|_2 \\
& & \mbox{}+\delta_2\|T_1^{-1}( (A^TW^TY)\otimes(A^TYW)
)\|_2+\delta_3\|T_1^{-1}\|_2.
\end{eqnarray*}

By setting $\delta_1=\|A\|_F$, $\delta_2=\|Q\|_F$ and
$\delta_3=\|G\|_F$, we can prove that
$$
\kappa_U(\psi^\Re) \leq \|\cals_1\|_2 \kappa_1^U(\psi^\Re).
$$

\section{\bf Structured Mixed and Componentwise
Condition Numbers} \label{seccomponentnorm}

Componentwise analysis \cite{CDW07,Higham94,Rohn89,Skeel79}
is more informative than its normwise
counterpart when the data are badly scaled or sparse.
Here, we consider the two
kinds of condition numbers introduced by Gohberg and Koltracht
\cite{GK93}. The first kind, called the mixed
condition number, measures the output errors in norm while the
input perturbations componentwise. The second kind, called the
componentwise condition number, measures both the output and the
input perturbations componentwise.

Following \cite{GK93}, the structured mixed and componentwise  condition
numbers for the complex CARE \eqref{care} are defined by
\begin{eqnarray*}
m (\varphi) &=& \lim_{\epsilon \to 0}
\sup_{\substack{\Delta \in C_\epsilon \\
\Delta A \in \C^{n\times n},\
\Delta G, \, \Delta Q \in \HM^{n\times n}} \atop
{G + \Delta G, \, Q + \Delta Q \succeq 0}}
\frac{\|\Delta X\|_{\max}}{\epsilon \|X\|_{\max}}, \\
c (\varphi) &=& \lim_{\epsilon \rightarrow 0}
\sup_{\substack{\Delta \in C_\epsilon \\
\Delta A \in \C^{n\times n},\
\Delta G, \, \Delta Q \in \HM^{n\times n}} \atop
{G + \Delta G, \, Q + \Delta Q \succeq 0}}
\frac{1}{\epsilon}\left\| \Delta X \oslash X \right\|_{\max} ,
\end{eqnarray*}
where $\Delta $ is defined in \eqref{eq:Delta all} and
\begin{align}
C_\epsilon &= \big\{\Delta A,\, \Delta G,\, \Delta Q ~|~| \Re(\Delta A))| \leq \epsilon |\Re(A)|,
|\Im(\Delta A))|\leq \epsilon |\Im(A)|, \nonumber \\
&  \quad |\sym(\Re(\Delta G))| \leq \epsilon |\sym(\Re(G))|,  |\ske(\Im(\Delta G))| \leq \epsilon |\ske(\Im(G))|, \nonumber\\
&  \quad |\sym(\Re(\Delta Q))| \leq \epsilon |\sym(\Re(Q))|, |\ske(\Im(\Delta Q))| \leq \epsilon |\ske(\Im(Q))|
\big\}. \label{eq:c eps}
\end{align}

In Theorems \ref{mcthm} and  \ref{th:4}, we present the structured mixed and
componentwise condition numbers of the complex CARE (\ref{care})
and the complex DARE (\ref{dare}).

\begin{theorem}\label{mcthm}
For the structured mixed and componentwise condition numbers of
the complex CARE (\ref{care}), we have respectively
\begin{align*}
&  m(\varphi) = \|X\|_{\max}^{-1}
\Big\| |Z^{-1}((I_n\otimes X) +(X^T\otimes I_n)\Pi)| \vect(|\Re(A)|) \\
&  \mbox{} + |Z^{-1}((I_n\otimes X) -(X^T\otimes I_n)\Pi)| \vect(|\Im(A)|)  \mbox{} + |Z^{-1}(X^T\otimes X)\cals_1|\sym(|\Re(G)|)\\
&  \mbox{} + |Z^{-1}(X^T\otimes X)\cals_2|\ske(|\Im(G)|)  \mbox{} + |Z^{-1}\cals_1|\sym(|\Re(Q)|) \mbox{} + |Z^{-1}\cals_2|\ske(|\Im(Q)|) \Big\|_\infty ,\\
&  c(\varphi) = \Big\| \diag(\vect(X))^\dagger  \cdot \Big( |Z^{-1}((I_n\otimes X) +(X^T\otimes I_n)\Pi)|
\vect(|\Re(A)|) \\
&  \mbox{} +|Z^{-1}((I_n\otimes X) -(X^T\otimes I_n)\Pi)| \vect(|\Im(A)|)  \mbox{} +|Z^{-1}(X^T\otimes X)\cals_1|\sym(|\Re(G)|) \\
&  \mbox{} +|Z^{-1}(X^T\otimes X)\cals_2|\ske(|\Im(G)|) \mbox{} +|Z^{-1}\cals_1|\sym(|\Re(Q)|)  \mbox{} +|Z^{-1}\cals_2|\ske(|\Im(Q)|) \Big) \Big\|_\infty,
\end{align*}
where $A^\dagger$ is the Moore-Penrose inverse of $A$. Furthermore, we have their simpler upper bounds
\begin{align*}
&m_U(\varphi):= \|X\|_{\max}^{-1} \|Z^{-1}\|_{\infty}
\Big\|\,|X|\, |\Re(A)|+|\Re(A)|^T |X| \\
& \mbox{} +|X|\, |\Im(A)|+|\Im(A)|^T |X| + |X| |\Re(G)||X| \mbox{} +|X|\, |\Im(G)|\,|X|+|\Re(Q)|+|\Im(Q)|\,\Big\|_{\max}
\end{align*}
and
\begin{align*}
&c_U(\varphi)
:= \|\diag(\vect(X))^\dagger Z^{-1}\|_{\infty} \cdot \Big\|\,|X|\, |\Re(A)|+|\Re(A)|^T |X| \\
&  \mbox{}+|X|\, |\Im(A)|+|\Im(A)|^T |X|+|X|\, |\Re(G)||X| \mbox{}+|X|\, |\Im(G)|\, |X|+|\Re(Q)|+|\Im(Q)|\,\Big\|_{\max}.
\end{align*}
\end{theorem}
{\bf Proof.}
From Lemma~\ref{lemma:care}, dropping the second and higher order terms,
we have
\begin{eqnarray*}
\|\Delta X\|_{\max} &=&
\|\vect(\Delta X)\|_\infty \approx \|Z^{-1}\calm \Delta \|_\infty\leq\|Z^{-1}\calm D_m\|_\infty \| D_m^\dagger \Delta \|_\infty,
\end{eqnarray*}
where
\begin{eqnarray*}
D_m &=& \diag ([\vect(\Re(A))^T,\ \vect(\Im(A))^T,  \, \sym(\Re(G))^T,\ \ske(\Im( G))^T, \\
&&\qquad \sym(\Re(Q))^T,\ \ske(\Im( Q))^T]^T) .
\end{eqnarray*}
Since $\Delta \in C_\epsilon$, we have
$\| D_m^\dagger \Delta \|_\infty \leq \epsilon$.
Because $\Delta$ can be chosen arbitrarily, the upper bound is
attainable. Recalling that $\bfe$ is the vector consisting
all 1's, it can be verified that
\begin{eqnarray*}
& & \left\|Z^{-1}\calm D_m\right\|_\infty=
\left\| |Z^{-1}\calm | \cdot |D_m|\bfe \right\|_\infty \\
&=& \Big\| |Z^{-1}\calm |\cdot
\Big[ |\vect(\Re(A))|^T,|\vect(\Im(A))|^T,|\sym(\Re(G))|^T, |\ske(\Im( G))|^T, \\
& & \qquad |\sym(\Re(Q))|^T, |\ske(\Im( Q))|^T\Big]^T\Big\|_\infty.
\end{eqnarray*}
After some algebraic manipulation, we can get the explicit
expression for $m(\varphi)$. From the following inequalities,
\begin{align*}
& |Z^{-1}((I_n\otimes X) +(X^T\otimes I_n)\Pi)|\vect(|\Re(A)|)\\
&  \leq |Z^{-1}|((I_n\otimes |X|) +(|X^T| \otimes I_n)\Pi)\vect(|\Re(A)|) =|Z^{-1}|\vect(|X| |\Re(A)|+|\Re(A)|^T |X| ),\\
&  |Z^{-1}((I_n\otimes X) -(X^T\otimes I_n)\Pi)| \vect(|\Im(A)|)
\leq|Z^{-1}|\vect(|X| |\Im(A)|+|\Im(A)|^T |X| ),\\
&  |Z^{-1}(X^T\otimes X)\cals_1|\sym(|\Re(G)|)\\
&\leq|Z^{-1}|(|X^T| \otimes |X|) \vect(|\Re(G)|) =|Z^{-1}| \vect(|X| |\Re(G)||X|),\\
& |Z^{-1}(X^T\otimes X)\cals_2|\ske(|\Im(G)|)\\
&\leq|Z^{-1}|(|X^T|\otimes |X|)\vect(|\Im(G)|) = |Z^{-1}| \vect(|X| |\Im(G)||X|),\\
&  |Z^{-1}\cals_1|\sym(|\Re(Q)|) \leq  |Z^{-1}| \vect(|\Re(Q)|),\quad  |Z^{-1}\cals_2|\ske(|\Im(Q)|) \leq  |Z^{-1}|\vect(|\Im(Q)|),
\end{align*}
and the monotonicity of the infinity norm,
we can obtain the upper bound $m_U(\varphi)$.
For the structured componentwise condition number $c(\varphi)$,
noting that
$$
\left\| \Delta X \oslash X \right\|_{\max} =
\| \diag( \vect(X))^\dagger \vect(\Delta X) \|_{\infty},
$$
similarly to the derivation of $m(\varphi)$, we can obtain
the explicit expression for $c(\varphi)$.
Also, the upper bound $c_U(\varphi)$ can be deduced similarly. \qed

The above theorem shows that an ill-conditioned $Z$
is an indication of large $m(\varphi)$ or $c(\varphi)$.

For the real CARE (\ref{Rcare}), the structured mixed and
componentwise condition
numbers of $\varphi^\Re$ at $X$ can be defined by
\begin{eqnarray*}
m (\varphi^\Re) &=& \lim_{\epsilon \to 0}
\sup_{\substack{\Delta_1 \in C_\epsilon^1 \\
\Delta A \in \C^{n\times n},\
\Delta G, \, \Delta Q \in \SM^{n\times n}} \atop
{G + \Delta G, \, Q + \Delta Q \succeq 0}}
\frac{\|\Delta X\|_{\max}}{\epsilon \|X\|_{\max}},\\
c (\varphi^\Re) &=& \lim_{\epsilon \rightarrow 0}
\sup_{\substack{\Delta_1 \in C_\epsilon^1 \\
\Delta A \in \C^{n\times n},\
\Delta G, \, \Delta Q \in \SM^{n\times n}} \atop
{G + \Delta G,\, Q + \Delta Q \succeq 0}}
\frac{1}{\epsilon}
\left\| \Delta X \oslash X \right\|_{\max},
\end{eqnarray*}
where $\Delta_1=[\Delta A^T, \sym(\Delta G)^T, \sym(\Delta Q)^T]^T$
and
\begin{eqnarray} \label{eq:c1 eps}
C_\epsilon^1 &=& \big\{
\Delta A,\, \Delta G, \, \Delta Q, ~|~ |\Delta A|\leq \epsilon |A|, \,
|\sym(\Delta G)| \leq  \epsilon |\sym(\Delta G)|, \nonumber \\
& & \quad |\sym(\Delta Q)|\leq \epsilon |\sym(\Delta Q)|\big\}.
\label{eq:C1 eps}
\end{eqnarray}

Now, we have expressions for $m(\varphi^\Re)$ and
$c(\varphi^\Re)$.
\begin{corollary}\label{mcthmRcare}
For structured mixed and componentwise condition numbers of the real CARE
(\ref{Rcare}), we have respectively
\begin{eqnarray*}
& & m(\varphi^\Re) = \|X\|_{\max}^{-1}
\Big\|| Z_1^{-1}((I_n\otimes X) +(X^T\otimes I_n)\Pi)| \vect(|A|) \\
& & \mbox{} + |Z_1^{-1}(X^T\otimes X)\cals_1|\sym(|G|) +
|Z_1^{-1}\cals_1|\sym(|Q|)\Big\|_\infty\\
& & c(\varphi^\Re) = \Big\|\diag(\vect(X))^\dagger
\Big(|Z_1^{-1}((I_n\otimes X) +(X^T\otimes I_n)\Pi)|\\
& & \cdot \vect(|A|) + |Z_1^{-1}(X^T\otimes X)\cals_1|\sym(|G|) \mbox{} +|Z_1^{-1}\cals_1|\sym(|Q|)\Big)\Big\|_\infty,
\end{eqnarray*}
where $Z_1=I_n\otimes (A-GX)^T+(A-GX)^T \otimes I_n$. Furthermore, we have their simpler upper bounds
\begin{eqnarray*}
m_U(\varphi^\Re) &:=& \|X\|_{\max}^{-1} \|Z_1^{-1}\|_{\infty}
\Big\|\,|X|\, |A|+|A|^T |X|\mbox{}+|X|\, |G|\, |X| + |Q|\,\Big\|_{\max},
\end{eqnarray*}
and
\begin{eqnarray*}
c_U(\varphi^\Re) &:=& \|\diag(\vect(X))^\dagger Z_1^{-1}\|_{\infty}
\Big\|\,|X|\, |A|+|A|^T |X|\mbox{} + |X|\, |G|\, |X| + |Q|\,\Big\|_{\max}.
\end{eqnarray*}
\end{corollary}

Zhou et al. \cite{zhou} derived the mixed and componentwise
condition numbers for the real CARE without exploiting the
symmetry structure. It can be readily
shown that our structured mixed and componentwise condition numbers
are smaller than their counterparts in \cite{zhou}, although,
they are empirically comparable.

Similarly, the structured mixed and componentwise  condition
numbers for the complex DARE \eqref{dare} can be defined by
\begin{eqnarray*}
m (\psi) &=& \lim_{\epsilon \to 0}
\sup_{\substack{\Delta \in C_\epsilon \\
\Delta A \in \C^{n\times n},\
\Delta G, \, \Delta Q \in \HM^{n\times n}} \atop
{G + \Delta G,\, Q + \Delta Q \succeq 0}}
\frac{\|\Delta Y\|_{\max}}{\epsilon \|Y\|_{\max}}, \\
c (\psi) &=& \lim_{\epsilon \rightarrow 0}
\sup_{\substack{\Delta \in C_\epsilon \\
\Delta A \in \C^{n\times n},\
\Delta G, \, \Delta Q \in \HM^{n\times n}} \atop
{G + \Delta G,\, Q + \Delta Q \succeq 0}}
\frac{1}{\epsilon}\left\| \Delta Y \oslash Y \right\|_{\max},
\end{eqnarray*}
where $\Delta $ is defined in \eqref{eq:Delta all} and
$C_\epsilon$ is defined in \eqref{eq:c eps}.

Following the proof of Theorem~\ref{mcthm}, we have
the structured mixed and componentwise condition numbers
for the complex DARE (\ref{dare}) in the following theorem. Its proof is omitted.

\begin{theorem}\label{th:4}
	With the notations as before, we have \begin{eqnarray*}
& & m(\psi)
= \|Y\|_{\max}^{-1}\Big\| |T^{-1}(I_n\otimes(A^H YW)  \mbox{}+((A^TW^TY^T)\otimes I_n)\Pi)|\cdot \vect(|\Re(A)|)\\
& & \mbox{}+|T^{-1}(I_n\otimes(A^H YW)-((A^TW^TY^T)\otimes I_n)\Pi)| \vect(|\Im(A)|)\\
& & \mbox{}+|T^{-1}(A^TW^TY^T)\otimes(A^HYW)\cals_1|\sym(|\Re(G)|)\\
& & \mbox{}+|T^{-1}(A^TW^TY^T)\otimes(A^HYW)\cals_2|\ske(|\Im(G)|)\\
& & \mbox{}+|T^{-1}\cals_1|\sym(|\Re(Q)|)+
|T^{-1}\cals_2|\ske(|\Im(Q)|)\Big\|_\infty
\end{eqnarray*}
and
\begin{eqnarray*}
& &c(\psi) =\Big\|\diag(\vect(Y))^\dagger \Big(|T^{-1}((I_n\otimes(A^H YW))  \mbox{}+((A^TW^TY^T)\otimes I_n)\Pi)|\cdot \vect(|\Re(A)|)\\
& & \mbox{}+|T^{-1}((I_n\otimes(A^H YW)) -((A^TW^TY^T)\otimes I_n)\Pi)| \vect(|\Im(A)|)\\
& & \mbox{}+|T^{-1}((A^TW^TY^T)\otimes(A^HYW)\cals_1 ) |\sym(|\Re(G)|)\\
& & \mbox{}+|T^{-1}((A^TW^TY^T)\otimes(A^HYW)\cals_2 ) |\ske(|\Im(G)|)\\
& & \mbox{}+|T^{-1}\cals_1|\sym(|\Re(Q)|)+
|T^{-1}\cals_2|\ske(|\Im(Q)|)\Big)\Big\|_\infty .
\end{eqnarray*}
Furthermore, we have their simpler upper bounds:
\begin{eqnarray*}
& & m_U(\psi) \\
&=& \|Y\|_{\max}^{-1} \|T^{-1}\|_{\infty}\Big\|\,|A^H||Y| |W| |\Re(A)|\\
& & \mbox{}+|\Re(A)|^T |Y||W||A|+|A^H||Y| |W| |\Im(A)|\\
& & \mbox{}+|\Im(A)|^T |Y||W||A|+|A^H||Y| |W| |\Re(G)||Y| |W| |A|\\
& & \mbox{}+|A^H||Y| |W| |\Im(G)||Y| |W| |A|+|\Re(Q)|+|\Im(Q)|\,\Big\|_{\max}
\end{eqnarray*}
and
\begin{eqnarray*}
& & c_U(\psi) \\
&=& \|\diag(\vect(Y))^\dagger T^{-1}\|_{\infty}
\Big\|\,|A^H||Y| |W| |\Re(A)|\\
& & \mbox{}+|\Re(A)|^T |Y||W||A|+|A^H||Y| |W| |\Im(A)|\\
& & \mbox{}+|\Im(A)|^T |Y||W||A|+|A^H||Y| |W| |\Re(G)||Y| |W| |A|\\
& & \mbox{}+|A^H||Y| |W| |\Im(G)||Y| |W| |A| +|\Re(Q)|+|\Im(Q)|\,\Big\|_{\max}.
\end{eqnarray*}

\end{theorem}
The above expressions in Theorem \ref{th:4} show that an ill conditioned
$T$ indicates large $m(\psi)$ or $c(\psi)$.

For the real DARE (\ref{Rdare}), we can define
the structured mixed and componentwise condition
numbers of $\psi^\Re$ at $Y$:
\begin{eqnarray*}
m (\psi^\Re) &=& \lim_{\epsilon \to 0}
\sup_{\substack{\Delta_1 \in C_\epsilon^1 \\
\Delta A \in \C^{n\times n},\
\Delta G,\, \Delta Q \in \SM^{n\times n}} \atop
{G + \Delta G,\, Q + \Delta Q \succeq 0}}
\frac{\|\Delta Y\|_{\max}}{\epsilon \|Y\|_{\max}}, \\
c (\psi^\Re) &=& \lim_{\epsilon \to 0}
\sup_{\substack{\Delta_1 \in C_\epsilon^1 \\
\Delta A \in \C^{n\times n},\
\Delta G,\, \Delta Q \in \SM^{n\times n}} \atop
{Q+\Delta Q \succeq 0,\, G + \Delta G,\, Q + \Delta Q \succeq 0}}
\frac{1}{\epsilon}\left\| \Delta Y \oslash Y \right\|_{\max},
\end{eqnarray*}
respectively,
where $\Delta_1=[\Delta A^T, \sym(\Delta G)^T, \sym(\Delta Q)^T]^T$
and $C_\epsilon^1$ is defined in \eqref{eq:c1 eps}.

Similarly to Corollary \ref{mcthmRcare}, we then can obtain the
structured mixed and componentwise condition numbers of the real DARE in the following corollary.
\begin{corollary}
With the notations above, we have
	\begin{eqnarray*}
m(\psi^\Re) &=&
\|Y\|_{\max}^{-1}\Big\||T^{-1}((I_n\otimes(A^T YW)) +((A^TYW)\otimes I_n)\Pi)|\cdot \vect(|A|)\\
& & \mbox{}+|T^{-1}((A^TYW)\otimes(A^TYW)\cals_1 ) |\sym(|G|)+|T^{-1}\cals_1|\sym(|Q|)\Big\|_\infty , \\
c(\psi^\Re) &=&
\Big\|\diag(\vect(Y))^\dagger \Big(|T^{-1}((I_n\otimes(A^T YW)) +((A^TYW)\otimes I_n)\Pi)|\cdot \vect(|A|)\\
& & \mbox{}+|T^{-1}((A^TYW)\otimes(A^TYW)\cals_1 ) |\sym(|G|)+|T^{-1}\cals_1|\sym(|Q|)\Big)\Big\|_\infty,
\end{eqnarray*}
Furthermore, we have their simpler upper bounds:
\begin{eqnarray*}
m_U(\psi^\Re)
&=& \|Y\|_{\max}^{-1} \|T^{-1}\|_{\infty}\Big\|\,|A^T||Y| |W| |A|+|A|^T |Y||W||A|\\
&&+|A^T||Y| |W| |G||Y| |W| |A|  +|Q|\,\Big\|_{\max}, \\
c_U(\psi^\Re)
&=& \|\diag(\vect(Y))^\dagger T^{-1}\|_{\infty}
\Big\|\,|A^T||Y| |W| |A| +|A|^T |Y||W||A|\\
&&+|A^T||Y| |W| |G||Y| |W| |A| +|Q|\,\Big\|_{\max}.
\end{eqnarray*}
\end{corollary}

\section{Small Sample Condition Estimation} \label{secsce}

Although the expressions of the condition numbers presented
earlier are explicit, they involve the solution and their
computation is intensive when the problem size is large.
Thus, practical algorithms for approximating the condition
numbers are worth studying \cite[Page 260]{Sun02}.
In this section, based on a small sample statistical condition
estimation method, we present a practical method for estimating
the condition numbers for the symmetric algebraic Riccati equations.

We first briefly describe our method.
Given a differentiable function
$f:\R^{p}\rightarrow \R$, we are interested in its sensitivity
at some input vector $x$. From its Taylor expansion, we have
$$
f(x+ \delta d)-f(x)= \delta (\nabla f(x))^T d+O(\delta ^2),
$$
for a small scalar $\delta$, where
\[
\nabla f(x)= \left[ \frac{\partial f(x)}{\partial x_1},
\frac{\partial f(x)}{\partial x_2},\ldots,
\frac{\partial f(x)}{\partial x_p} \right]^T
\]
is the gradient of $f$ at $x$. Then the local sensitivity, up to the first order
in $\delta$, can be measured by
$\|\nabla f(x)\|_2$. The condition number of $f$ at $x$ is mainly
determined by the norm of the gradient $\nabla f(x)$
(\cite{KenneyLaub_SISC94}). It is shown in \cite{KenneyLaub_SISC94}
that if we select $d$ uniformly and randomly from the unit $p$-sphere
$S_{p-1}$ (denoted $d \in U(S_{p-1})$), then the expected value
${\bf E}(|(\nabla f(x))^T d|/\omega_p)$ is $\|\nabla f(x)\|_2$,
where $\omega_p$ is the Wallis factor,
which depends only on $p$, given by
$$
\omega_p=
\begin{cases}
1, & \text{for}~p\equiv 1,\\
\frac{2}{\pi}, &\text{for}~p\equiv 2,\\
\frac{1 \cdot 3 \cdot 5 \cdots (p-2)}{2 \cdot 4 \cdot 6 \cdots (p-1)}, &
\text{for}~p ~ \text{odd} ~\text{and} ~p>2, \\
\frac{2}{\pi} \frac{2 \cdot 4 \cdot 6 \cdots (p-2)}
{1 \cdot 3 \cdot 5 \cdots (p-1)}, &
\text{for}~ p ~ \text{even}~ \text{and} ~p>2 ,
\end{cases}
$$
which can be accurately
approximated by
\begin{equation}\label{ss:p} \omega_p\approx
\sqrt{\frac{2}{\pi(p-\frac{1}{2})}}.
\end{equation}
Therefore,
$$
\nu=\frac{|(\nabla f(x))^T d|}{\omega_p}
$$
can be used to estimate $\|\nabla f(x)\|_2$, an approximation
of the condition number, with high probability
\cite{KenneyLaub_SISC94}. Specifically, for $\gamma >1$,
\begin{eqnarray*}
 \Prob\left(\frac{\|\nabla f(x)\|_2}{\gamma}\leq \nu \leq
\gamma \|\nabla f(x)\|_2\right)
\geq 1-\frac{2}{\pi \gamma}+O(\gamma^{-2}).
\end{eqnarray*}
Multiple samples $d_j$ can be used to increase the accuracy
\cite{KenneyLaub_SISC94}. The $k$-sample condition estimation
is given by
\begin{eqnarray*}
 \nu(k)
= \frac{\omega_k}{\omega_p}
\sqrt{|\nabla f(x)^T d_1|^2 +|\nabla f(x)^T d_2|^2 + \cdots +
|\nabla f(x)^T d_k|^2},
\end{eqnarray*}
where $d_1,d_2,...,d_k$ are orthonormalized after they
are selected uniformly and randomly from $U(S_{p-1})$.
In particular, the accuracy of $\nu(2)$ is given by
$$
\Prob\left(\frac{\|\nabla f(x)\|_2}{\gamma}\leq \nu(2) \leq
\gamma \|\nabla f(x)\|_2\right)
\approx 1-\frac{\pi}{4 \gamma^2}.
$$
Usually, a small set of samples is sufficient for good accuracy.

These results can be readily generalized
to vector-valued or matrix-valued functions by viewing $f$ as
a map from $\R^s$ to $\R^t$ by applying the operations
$\vect$ and $\sf unvec$ to transform data between matrices and vectors,
where each of the $t$ entries of $f$ is a scalar-valued function.

\subsection{Structured normwise case}

In this subsection, by applying the small sample condition
estimation method, we devise the algorithms for estimating
the structured normwise condition numbers
of CARE and DARE. Before that, we introduce the {\sf
unvec} operation: Given $m$ and $n$, for
$v=[v_1,v_2,...,v_{mn}] \in \R^{1 \times mn}$,
$A={{\sf unvec}}(v)$ sets the $(i,j)$-entry of $A$ to
$v_{i+(j-1)n}$.

For CARE, from Lemma~\ref{lemma:care}, the directional derivative
$D_X \in \HM^{n\times n}$ of $\varphi$ at $X$ with respect to
the direction $\Delta$, defined in \eqref{eq:Delta all},
satisfies the continuous Lyapunov equation (\ref{eq:Lyapunov}).
Putting things together,
we propose the following subspace structured condition number
estimation algorithm for the complex CARE \eqref{care}.

\begin{enumerate}
\item
Generate vectors $f_i \in \R^{4n^2}$, $i=1,...,k$,
with each entry in ${\mathcal N}(0,1)$. Orthonormalize them
using, for example, the QR factorization, to get $z_j\in \R^{4n^2}$,
$j=1,...,k$. Each $z_j$ can be converted into
matrices $\widetilde{A_j}$, $\widetilde{G_j}$,
and $\widetilde{Q_j}$ by applying the {\sf unvec}
operation, where $\widetilde{A_j} \in \C^{n\times n}$
and $\widetilde{G_j}, \, \widetilde{Q_j} \in \HM^{n\times n}$;
\item
For $i=1,2,...,k$, solve for $D_i \in \HM^{n\times n}$ in
the following continuous Lyapunov equation
\begin{eqnarray*}
(A-GX)^H D_i+D_i(A-GX)
= X\widetilde{G_i} X-\widetilde{Q_i}
-X\widetilde{A_i}-\widetilde{A_i}^H X;
\end{eqnarray*}
\item
Approximate $\omega_k$ and $\omega_p$ ($p=4n^2$)
by~(\ref{ss:p}) and calculate the absolute condition
number matrix
\begin{align*}
K_{\rm abs}^{\mathrm{CARE},(k)} :=\left\|[A,G,Q]\right\|_F\frac{\omega_k}{\omega_p}
\sqrt{|D_1|^2+ |D_2|^2+\cdots+|D_k|^2},
\end{align*}
where the square operation is applied to each
entry of $|D_i|$, $i=1,2,...,k$ and the square root
is also applied componentwise;
\item
Finally, the relative condition number matrix
$$
K_{\rm rel}^{\mathrm{CARE},(k)} =
K_{\rm abs}^{\mathrm{CARE},(k)} \oslash X
$$
is obtained by componentwise division for nonzero
entries of $X$, leaving the entries of
$K_{\rm abs}^{\mathrm{CARE},(k)}$ corresponding to
the zero entries of $X$ unchanged.
\end{enumerate}

The real CARE \eqref{Rcare} is a special case.

Note that Step 2 in the above algorithm involves solving a
sequence of Lyapunov equations. When a Lyapunov equation
is ill-conditioned, the computed solution $D_i$ can be
inaccurate, consequently, the condition number for CARE computed
in the following Step 3 can be inaccurate. However, the
conditioning of the Lyapunov equation and that of CARE are
related in that the ill-conditioning of the continuous Lyapunov
equation implies the ill-conditioning of the original CARE,
because solving the Lyapunov equation is essentially equivalent
to finding $Z^{-1}$ in the condition number for CARE presented
in Theorem \ref{mcthm}.


For the complex DARE \eqref{dare}, from Lemma \ref{lemma:dare},
the directional derivative $D_Y \in \HM^{n\times n}$ of $\psi$
at $Y$ with respect to the direction $\Delta$ is the solution
of the discrete Lyapunov equation
\begin{eqnarray*}
D_Y-(WA)^H D_Y W A
= \Delta Q+(A^HYW)\Delta A+\Delta A^H (YWA)-(A^HYW)\Delta G (YWA) ,
\end{eqnarray*}
where $\Delta A\in \C^{n\times n}$ and
$\Delta G, \, \Delta Q \in \HM^{n\times n}$.

Similarly to the complex CARE case, we propose the following
algorithm for the complex DARE.
\begin{enumerate}
\item
Generate vectors $f_i \in \R^{4n^2}$, $i=1,...,k$
with each entry in ${\mathcal N}(0,1)$. Orthonormalize them
using, for example, the QR factorization, to get $z_j\in \R^{4n^2}$,
$j=1,...,k$. Each $z_j$ can be converted into the corresponding
matrices $\widetilde{A_j}$, $\widetilde{G_j}$,
and $\widetilde{Q_j}$ by applying the {\sf unvec}
operation, where $\widetilde{A_j} \in \C^{n\times n}$
and $\widetilde{G_j}, \, \widetilde{Q_j} \in \HM^{n\times n}$;
\item
For $i=1,2,...,k$, solve for $D_i \in \HM^{n\times n}$ in
the following discrete Lyapunov equation
\begin{eqnarray*}
D_i-(WA)^H D_i WA
=\widetilde{Q_i}+(A^HYW)\widetilde{A_i}+\widetilde{A_i}^H (YWA)(A^HYW)\widetilde{G_i}(YWA);
\end{eqnarray*}
\item
Approximate $\omega_k$ and $\omega_p$ ($p=4n^2$)
by~(\ref{ss:p}) and calculate the absolute condition
number matrix
\begin{align*}
K_{\rm abs}^{\mathrm{DARE},(k)} :=\|[A,G,Q]\|_F
\frac{\omega_k}{\omega_p}
\sqrt{|D_1|^2+ |D_2|^2+\cdots+|D_k|^2};
\end{align*}
\item
Finally, the relative condition number matrix
$$
K_{\rm rel}^{\mathrm{DARE},(k)} =
K_{\rm abs}^{\mathrm{DARE},(k)} \oslash Y .
$$
\end{enumerate}


\subsection{Structured componentwise case}

Componentwise condition number often leads to a more realistic
indication of the accuracy of a computed solution than
the normwise condition number.
The sensitivity effects of componentwise perturbations can be measured by
the SCE method \cite{KenneyLaub_SISC94}.
For a perturbation $\Delta A=[\Delta a_{ij}]$ on a matrix
$A=[a_{ij}] \in \R^{m \times n}$, it
is a componentwise perturbation, if
$$
|\Delta A| \leq \e |A| \quad {\rm or } \quad |\Delta a_{ij}|\leq \e |a_{ij}|.
$$
We can write $\Delta A=\delta\cdot {\mathcal A} \boxdot A$ with
$|\delta|\leq \e$ and the entries of $\mathcal A$ are
in the interval $[-1,1]$, where $\boxdot$ is a componentwise
multiplication. We propose the following algorithm for
a structured componentwise sensitivity estimate of the
solution $X$ of the complex CARE \eqref{care}.
\begin{enumerate}
\item
Generate vectors $f_i \in \R^{4n^2}$, $i=1,...,k$,
with each entry in ${\mathcal N}(0,1)$. Orthonormalize them
using, for example, the QR factorization, to get $z_j\in \R^{4n^2}$,
$j=1,...,k$. Each $z_j$ can be converted into the corresponding
matrices $\widetilde{A_j}$, $\widetilde{G_j}$,
and $\widetilde{Q_j}$ by applying the {\sf unvec}
operation, where $\widetilde{A_j} \in \C^{n\times n}$
and $\widetilde{G_j}, \, \widetilde{Q_j} \in \HM^{n\times n}$;
\item
For $j=1,2,...,k$, set $[\widetilde{A_j},\ \widetilde{G_j},
\ \widetilde{Q_j}]$ equal to the componentwise product of
$[A,\ G,\ Q]$ and $[\widetilde{A_j},\ \widetilde{G_j},
\ \widetilde{Q_j}]$;
\item
For $i=1,2,...,k$, solve for $D_i \in \HM^{n\times n}$ in
the following continuous Lyapunov equation
\begin{eqnarray*}
 (A-GX)^H D_i+D_i(A-GX)
=X\widetilde{G_i} X-\widetilde{Q_i}
-X\widetilde{A_i}-\widetilde{A_i}^H X;
\end{eqnarray*}
\item
Approximate $\omega_k$ and $\omega_p$ ($p=4n^2$)
by~(\ref{ss:p}) and calculate the absolute condition
number matrix
$$
C_{\rm abs}^{\mathrm{CARE},(k)} :=
\frac{\omega_k}{\omega_p}
\sqrt{|D_1|^2+ |D_2|^2+\cdots+|D_k|^2};
$$
\item
Finally, the relative condition number matrix
$$
C_{\rm rel}^{\mathrm{CARE},(k)} =
C_{\rm abs}^{\mathrm{CARE},(k)} \oslash X .
$$
\end{enumerate}

Analogously to the above complex case, we propose
the following algorithm for the complex DARE \eqref{dare}.
\begin{enumerate}
\item
Generate vectors $f_i \in \R^{4n^2}$, $i=1...,k$,
with each entry in ${\mathcal N}(0,1)$. Orthonormalize them
using, for example, the QR factorization, to get $z_j\in \R^{4n^2}$,
$j=1,...,k$. Each $z_j$ can be converted into the corresponding
matrices $\widetilde{A_j}$, $\widetilde{G_j}$,
and $\widetilde{Q_j}$ by applying the {\sf unvec}
operation, where $\widetilde{A_j} \in \C^{n\times n}$
and $\widetilde{G_j}, \, \widetilde{Q_j} \in \HM^{n\times n}$;
\item
For $j=1,2,...,k$, set $[\widetilde{A_j},\ \widetilde{G_j},
\ \widetilde{Q_j}]$ equal to the componentwise product of
$[A,\ G,\ Q]$ and $[\widetilde{A_j},\ \widetilde{G_j},
\ \widetilde{Q_j}]$;
\item
For $i=1,2,...,k$, solve for $D_i \in \HM^{n\times n}$ in
the following discrete Lyapunov equation
\begin{eqnarray*}
 D_i-(WA)^H D_i WA
= \widetilde{Q_i}+(A^HYW)\widetilde{A_i}+\widetilde{A_i}^H (YWA) -(A^HYW)\widetilde{G_i}(YWA);
\end{eqnarray*}
\item
Approximate $\omega_k$ and $\omega_p$ ($p=4n^2$)
by~(\ref{ss:p}) and calculate the absolute condition
number matrix
$$
C_{\rm abs}^{\mathrm{DARE},(k)} :=
\frac{\omega_k}{\omega_p}
\sqrt{|D_1|^2+ |D_2|^2+\cdots+|D_k|^2};
$$
\item
Finally, the relative condition number matrix
$$
C_{\rm rel}^{\mathrm{DARE},(k)} =
C_{\rm abs}^{\mathrm{DARE},(k)} \oslash Y .
$$
\end{enumerate}

\section{\bf Numerical Examples}\label{secnumex}

\begin{table*}[t]
\begin{center}
\caption{Comparison of the accurate relative changes in
the solution with the estimates obtained by our condition numbers,
where $\epsilon = 10^{-8}$.}
\label{tblNorm}

\begin{tabular}{|c|ccc|}
\hline
$\nu $ & $\|\Delta X\|_F/\|X\|_F$ & $\epsilon \, \kappa_U(\varphi^\Re)$
& $\epsilon \, \kappa_1^U (\varphi^\Re)$ \\
\hline \hline
$1$ & $ 3.0642\times 10^{-9}$ & $ 3.7258 \times 10^{-8}$ &
$4.0054 \times 10^{-8}$ \\
\hline
$10^6$ & $  7.0865\times 10^{-9}$ & $5.000\times 10^{-3}$ &
$5.000\times 10^{-3}$  \\
\hline
$10^{-6}$ & $  4.6983\times 10^{-9}$ & $5.0000\times 10^{3}$ &
$5.0000\times 10^{3}$ \\
\hline
\end{tabular}

\begin{tabular}{|c|cccc|}
\hline
$\nu $ & $\|\Delta X\|_{\max}/\|X\|_{\max}$ &
$\epsilon \, m(\varphi^\Re)$ & $\left\|\Delta X\oslash X\right\|_{\max}$ &
$\epsilon \, c(\varphi^\Re)$\\
\hline \hline
$1$ &
$ 6.1630\times 10^{-9}$ & $ 1.6667 \times 10^{-8}$ &
$7.7288 \times 10^{-9}$ &$ 1.6667 \times 10^{-8}$ \\
\hline
$10^6$ & $  7.0865\times 10^{-9}$ & $1.5000 \times 10^{-8}$ &
$   1.2161\times 10^{-8}$& $1.5000 \times 10^{-8}$ \\
\hline
$10^{-6}$  &
$  4.6983\times 10^{-9}$ & $ 2.0000 \times 10^{-8}$ &
$  7.5086\times 10^{-9}$ & $ 2.0000 \times 10^{-8}$ \\
\hline
\end{tabular}
\end{center}
\end{table*}

In this section, we adopt the examples in \cite{BennerC,BennerD,zhou} to illustrate
the effectiveness of our methods. All the experiments were performed
using  \textsc{Matlab} 8.1, with the machine epsilon $\mu \approx
2.2 \times 10^{-16}$.

Given $A\in \R^{n\times n}$ and $G,\,Q\in \SM^{n\times n}$,
we generated the perturbations on $A$, $G$ and $Q$ as follows:
$\Delta A= \epsilon (M_1 \boxdot A)$,
$\Delta G= \epsilon (M_2 \boxdot G)$,
and $\Delta Q= \epsilon (M_3 \boxdot Q)$, where
$\epsilon = 10^{-j}$ for some nonnegative integer $j$,
$\boxdot$ denotes the componentwise
multiplication of two matrices, and $M_1 \in \R^{n\times n}$, and
$M_2, M_3 \in \SM^{n\times n}$ whose entries are random numbers
uniformly distributed in the open interval $(-1,1)$.

{\bf Example 1}
Consider the CARE (\ref{care}) from \cite[Example 9]{BennerC} with
\[
A=\begin{bmatrix}
  0&\nu\cr0&0
\end{bmatrix},\quad Q=I_2,\quad G=BR^{-1}B^T,
\]
where
\[
B=\begin{bmatrix}0\cr 1\end{bmatrix},\ R=1.
\]
The pair $(A,G)$ is c-stabilizable and the pair $(A,Q)$ is
c-detectable. The exact solution is
$$
X=\begin{bmatrix}
  \frac{\sqrt{1+2\nu}}{\nu} & 1\cr 1& \sqrt{1+2\nu}
\end{bmatrix}.
$$
When $\nu$ is large or small, $\|X\|_F$ is approximately
$\sqrt{\nu} \, (\nu\geq 1)$ or $1/\sqrt{\nu} \, (0<\nu<1)$ respectively and CARE becomes
ill conditioned in terms of the normwise conditions
$\kappa_U(\varphi^\Re)$
and $\kappa_1^U (\varphi^\Re)$. However, as shown in Table \ref{tblNorm},
from the componentwise perturbation analysis,
$m(\varphi^\Re)$ and $c(\varphi^\Re)$ are always of $\Oh(1)$.

Let $\tilde{Q}=Q+\Delta Q,\tilde{A}=A+\Delta A,\tilde{G}=G+\Delta G$
be the coefficient matrices of the perturbed CARE (\ref{perturbequ}).
The perturbation size $\epsilon = 10^{-8}$. We
used the \textsc{Matlab} function \texttt{are} to compute the unique symmetric positive semidefinite solution
$\tilde{X}$ to the perturbed equation.
Let $\Delta X = \tilde{X} - X$.

For the bound $\kappa_U(\varphi^\Re)$, we set $\delta_1=\|A\|_F,\delta_2=\|\sym(Q)\|_2,\delta_3=\|\sym(G)\|_2$. For $\kappa_1^U (\varphi^\Re)$ in \cite{zhou}
we choose $\delta_1=\|A\|_F,\delta_2=\|Q\|_F,\delta_3=\|G\|_F$.
Table~\ref{tblNorm} compares the accurate relative changes
$\| \Delta X \|_F/\|X\|_F$,
$\| \Delta X \|_{\max}/\|X\|_{\max}$ and
$\left\| \Delta X \oslash X\right\|_{\max}$
obtained by MATLAB with the estimates obtained
by our condition numbers. Our normwise condition numbers
are consistent with those in \cite[page 9]{BennerC} for
$\nu = 1, 10^{6}, 10^{-6}$. Our mixed and componentwise
condition numbers, however, give accurate estimates for
the corresponding relative changes in the solution.

For the SCE algorithms,
we set the sample number $k=5$ and tested them for various values of $\nu$.
The results are shown as follows.
For $\nu=1$,
\begin{align*}
\Delta X \oslash X&= 10^{-8 }\,\begin{bmatrix}
   -0.4039 &  -0.7729\\
   -0.7729 &  -0.6163
\end{bmatrix},\cr
\epsilon \, K_{\rm rel}^{\mathrm{CARE},(5)}&= 10^{-8}\,\begin{bmatrix}
    6.5364 &   7.4764\\
    7.4764   & 3.8048
\end{bmatrix} ,\cr
\epsilon \, C_{\rm rel}^{\mathrm{CARE},(5)}&=10^{-8} \, \begin{bmatrix}
     0.7649  &  0.6111\\
    0.6111    &0.6727
\end{bmatrix} .
\end{align*}
For $\nu=10^{6}$,
\begin{align*}
\Delta X \oslash X&=10^{-7 }\,\begin{bmatrix}
    -0.1216 &  -0.0962\\
   -0.0962  & -0.0709
\end{bmatrix},\cr
\epsilon \,K_{\rm rel}^{\mathrm{CARE},(5)}&=  10^{-2} \begin{bmatrix}
    0.9807  &  1.3288\\
    1.3288  &  1.1426\end{bmatrix} ,\cr
\epsilon \,C_{\rm rel}^{\mathrm{CARE},(5)}&= 10^{-8}\begin{bmatrix}
       1.1150  &  0.8994\\
    0.8994   & 1.2533
\end{bmatrix} .
\end{align*}
For $\nu=10^{-6}$,
\begin{align*}
\Delta X \oslash X&= 10^{-8 }\,\begin{bmatrix}
  -0.4698  & -0.4207\\
   -0.4207  & -0.7509
\end{bmatrix},\cr
\epsilon \,K_{\rm rel}^{\mathrm{CARE},(5)}&= 10^{4}\, \begin{bmatrix}
  1.0725 &   1.0725\\
    1.0725  &  0.0000\end{bmatrix} ,\cr
\epsilon \,C_{\rm rel}^{\mathrm{CARE},(5)}&= 10^{-8}\,\begin{bmatrix}
     1.1728 &   0.4670\\
    0.4670  &  0.7666
\end{bmatrix} .
\end{align*}
As we can see, for this particular example, the componentwise condition
matrices $C_{\rm rel}^{\mathrm{CARE},(5)}$ for all values of $\nu$
can be used to accurately estimate
the changes in the solution. In contrast, the normwise condition matrix
$K_{\rm rel}^{\mathrm{CARE},(5)}$ can give good estimate only
when $\nu=1$ because the problem is well conditioned under the normwise
perturbation analysis in this case.

{\bf Example 2}
For DARE, we  adopt the following example from \cite{zhou}.
Consider the DARE (\ref{dare}) with
\[
Q=VQ_0V,\quad A=VA_0V,\quad G=VG_0V,
\]
where
\[
Q_0=\diag([10^m,1,10^{-m}]^T), \quad
A_0=\diag([0,10^{-m},1]^T),
\]
\[
G_0=\diag([10^{-m},10^{-m},10^{-m}]^T),
\]
and
\[
V=I-2vv^T/3, \quad v=[1,1,1]^T.
\]
Correspondingly, in the original DARE (\ref{dareorigin}),
$B = V$, $R = G_0^{-1}$, and $C = V \sqrt{Q_0} V$.  The
pair $(A,B)$ is d-stabilizable and the pair $(A,C)$ is d-detectable.
The unique symmetric positive semidefinite solution $Y$ to the DARE
(\ref{dare}) is given by $Y=VY_0V$, where
$Y_0=\diag([y_1,y_2,y_3]^T)$ with
$$
y_i=(a_i^2+q_ig_i-1+((a_i^2+q_ig_i-1)^2+4q_ig_i)^{1/2})/(2g_i),
$$
and $q_i$, $a_i$ and $g_i$ are respectively the diagonal elements
of $Q_0$, $A_0$ and $G_0$. The perturbation matrices
$\Delta A$, $\Delta G$ and $\Delta Q$ were generated as
described in the beginning of this section with $\epsilon=10^{-12}$. Let $\tilde{Q}=Q+\Delta
Q,\tilde{A}=A+\Delta A,\tilde{G}=G+\Delta G$ be the coefficient
matrices of the perturbed DARE (\ref{dare}). We used \textsc{Matlab}
function \texttt{dare} to compute the
unique symmetric positive semidefinite solution $\tilde{Y}$ of the
perturbed equation (\ref{perturbequdare}). Let $\Delta Y = \tilde{Y} - Y$.

For the bound $\kappa_U(\psi^\Re)$, we set
$\delta_1=\|A\|_F,\delta_2=\|\sym(Q)\|_2,\delta_3=\|\sym(G)\|_2$.
For $\kappa_1^U (\psi^\Re)$ in \cite{zhou} we choose $\delta_1=\|A\|_F,\delta_2=\|Q\|_F,\delta_3=\|G\|_F$. Table~\ref{tblNormdare} shows that our condition numbers give
reasonably good estimates for the changes in the solution.



\begin{table*}
\begin{center}
\caption{Comparison of the accurate relative changes in
the solution with the estimates obtained by our condition numbers,
where $\epsilon = 10^{-12}$.}
\label{tblNormdare}

\begin{tabular}{|c|ccc|}
\hline $m $ & $\|\Delta Y\|_F/\|Y\|_F$ & $\epsilon \, \kappa_U(\psi^\Re)$
& $\epsilon \, \kappa_1^U (\psi^\Re)$ \\
\hline \hline
$1$ & $1.2974\times 10^{-12}$ & $6.6183\times 10^{-12}$ &
$7.1051\times 10^{-12}$ \\
\hline
$5$ & $1.5931\times 10^{-8}$ & $5.0002\times 10^{-8}$ &
$5.2934\times 10^{-8}$ \\
\hline
$7$ & $1.0577\times 10^{-7}$ & $5.0000\times 10^{-6}$ &
$5.2932\times 10^{-6}$  \\
\hline
\end{tabular}
\begin{tabular}{|c|cccc|}
\hline $m $  & $\|\Delta Y\|_{\max}/\|Y\|_{\max}$ &
$\epsilon \, m(\psi^\Re)$ & $\left\|\Delta Y\oslash Y \right\|_{\max}$ &
$\epsilon \, c(\psi^\Re)$\\
\hline \hline
$1$ & $1.2945\times 10^{-12}$ & $4.8227\times 10^{-12}$ &
$1.8177\times 10^{-12}$& $1.1056\times 10^{-11}$\\
\hline
$5$  & $1.5931\times 10^{-8}$ & $3.9507\times 10^{-8}$ &
$6.3715\times 10^{-8}$ &$1.5801\times 10^{-7}$ \\
\hline
$7$  & $1.0577\times 10^{-7}$ & $3.9506\times 10^{-6}$ &
$4.2307\times 10^{-7}$ & $1.5802\times 10^{-5}$ \\
\hline
\end{tabular}
\end{center}
\end{table*}


For the SCE algorithms,
we set the sample number $k=5$ and tested them for various values of $m$ with $\epsilon=10^{-12}$.
The results are shown as follows.
For $m=1$,
\begin{align*}
\Delta Y \oslash Y&= 10^{-11 }\,\begin{bmatrix}
     -0.1465  &  0.1818 &  -0.1187\\
    0.1818 &  -0.1294  &  0.1432\\
   -0.1187  &  0.1432  & -0.0856
\end{bmatrix},\cr
\epsilon \,K_{\rm rel}^{\mathrm{DARE},(5)}&=10^{-12 }\, \begin{bmatrix}
 25.5340  &-25.0489 & -15.4023\\
  -25.0489 &  10.8400  &  6.1395\\
  -15.4023 &   6.1395  &  3.1528\end{bmatrix} ,\cr
\epsilon \,C_{\rm rel}^{\mathrm{DARE},(5)}&=10^{-12 }\, \begin{bmatrix}
    3.0920&   -3.1064  & -2.1995\\
   -3.1064  &  1.4491  &  0.3301\\
   -2.1995   & 0.3301   & 0.9229
\end{bmatrix} .
\end{align*}
For $m=5$,
\begin{align*}
\Delta Y \oslash Y&=10^{-7 }\, \begin{bmatrix}
   -0.6371 &   0.3186 &  -0.1593\\
    0.3186 &  -0.1593  &  0.0797\\
   -0.1593  &  0.0797 &  -0.0398
\end{bmatrix},\cr
\epsilon \,K_{\rm rel}^{\mathrm{DARE},(5)}&= 10^{-7}\, \begin{bmatrix}
  5.1319  & -2.5662 &  -1.2831\\
   -2.5662  &  1.2831 &   0.6415\\
   -1.2831   & 0.6415 &   0.3208
\end{bmatrix} ,\cr
\epsilon \,C_{\rm rel}^{\mathrm{DARE},(5)}&=  10^{-8}\, \begin{bmatrix}
      6.1758 &  -3.0882 &  -1.5441\\
   -3.0882  &  1.5441  &  0.7721\\
   -1.5441  &  0.7721  &  0.3860
\end{bmatrix} .
\end{align*}
For $m=7$,
\begin{align*}
\Delta Y \oslash Y&=10^{-6 }\,\begin{bmatrix}
     0.4231 &  -0.2115  &  0.1058\\
   -0.2115  &  0.1058  & -0.0529\\
    0.1058  & -0.0529 &   0.0264
\end{bmatrix},\cr
\epsilon \,K_{\rm rel}^{\mathrm{DARE},(5)}&= 10^{-5} \, \begin{bmatrix}
 4.0438 &  -2.0219  & -1.0109\\
   -2.0219  &  1.0109  &  0.5055\\
   -1.0109  &  0.5055  &  0.2527\end{bmatrix} ,\cr
\epsilon \,C_{\rm rel}^{\mathrm{DARE},(5)}&= 10^{-5} \, \begin{bmatrix}
      1.0597 &  -0.5299 &  -0.2649\\
   -0.5299 &   0.2649  &  0.1325\\
   -0.2649 &   0.1325  &  0.0662
\end{bmatrix} .
\end{align*}

As shown above, even for a small number of samples,
the accuracy of the SCE method is within
a factor between $10^{-1}$ and $10$, which is considered
acceptable \cite[Chapter 15]{Higham02}.

\section{\bf Concluding Remarks}\label{secconclude}

In this paper, by exploiting the symmetry structure
and separating the real and imaginary parts,
we present structured perturbation
analyses of both the continuous-time and the discrete-time symmetric
algebraic Riccati equations. From the analyses, we define
the structured
normwise, mixed and componentwise condition numbers and derive
their upper bounds. Our bounds are improvements of the
results in previous work \cite{Sun02,zhou}.
Our preliminary experiments show that the three kinds of condition
numbers provide accurate bounds for the change in the perturbed solution.
Also, applying the small-sample
condition estimation method, we propose statistical
algorithms for practically estimating the structured condition numbers
for continuous and discrete symmetric algebraic Riccati equations.

\section*{Acknowledgments}

The authors would like to thank the reviewer for his/her comments, which  improve the presentation of the earlier version of this paper. S. Qiao is partially supported by
Natural Science and Engineering Council (NSERC) of Canada.


\end{document}